\documentclass[a4paper,11pt]{article}

\pdfoutput=1 

\usepackage[margin=0.75in]{geometry}

\usepackage[utf8]{inputenc}
\usepackage[T1]{fontenc}
\usepackage{lmodern}

\usepackage{amsmath,amssymb,amsthm,mathrsfs}
\usepackage{hyperref}
\usepackage[capitalise,noabbrev]{cleveref}
\usepackage{microtype}

\usepackage{listings}
\usepackage{xcolor}

\definecolor{keywordcolor}{RGB}{0,0,180}
\definecolor{commentcolor}{RGB}{128,128,128}
\definecolor{stringcolor}{RGB}{0,128,0}

\lstdefinelanguage{lean4}{
  keywords={def, theorem, lemma, example, structure, class, instance, where, if, then, else, let, in, do, return, match, with, fun, forall, exists, Prop, Type, Sort, noncomputable, variable, section, namespace, end, open, import, axiom, constant, inductive, extends, deriving},
  sensitive=true,
  morecomment=[l]{--},
  morecomment=[s]{/-}{-/},
  morestring=[b]",
  literate={α}{{$\alpha$}}1 {β}{{$\beta$}}1 {γ}{{$\gamma$}}1 {δ}{{$\delta$}}1 {ε}{{$\varepsilon$}}1
           {τ}{{$\tau$}}1 {φ}{{$\varphi$}}1 {ℝ}{{$\mathbb{R}$}}1 {ℕ}{{$\mathbb{N}$}}1
           {→}{{$\to$}}1 {←}{{$\leftarrow$}}1 {↦}{{$\mapsto$}}1 {∈}{{$\in$}}1 {∀}{{$\forall$}}1
           {∃}{{$\exists$}}1 {∧}{{$\land$}}1 {∨}{{$\lor$}}1 {¬}{{$\lnot$}}1 {≤}{{$\leq$}}1
           {≥}{{$\geq$}}1 {≠}{{$\neq$}}1 {·}{{$\cdot$}}1 {×}{{$\times$}}1 {‖}{{$\|$}}1
           {₀}{{$_0$}}1 {₁}{{$_1$}}1 {₂}{{$_2$}}1 {ₓ}{{$_x$}}1
           {𝓘}{{$\mathscr{I}$}}1 {∫}{{$\int$}}1 {≥}{{$\geq$}}1
}

\lstnewenvironment{leancode}{\lstset{language=lean4}}{}
\newcommand{\lean}[1]{\lstinline[language=lean4]{#1}}

\theoremstyle{definition}

\theoremstyle{remark}

\title{Integral Curves and Flows on Banach Manifolds in Lean}

\author{
  Weichen Winston Yin\thanks{University of California, Berkeley. Email: winstonyin@berkeley.edu. ORCID: 0000-0001-5692-7671}
  \and
  Yury Kudryashov\thanks{Harmonic. Email: urkud@urkud.name. ORCID: 0000-0003-4286-9276}
}

\date{}

\begin{document}

\maketitle

\begin{abstract}
  We present a formalisation of the existence and uniqueness theorems of integral curves of vector fields on Banach manifolds in the Lean theorem prover. First, we formalize properties of differential equations on Banach spaces (the Picard-Lindel\"of theorem, the Gr\"onwall inequality, and corollaries), and then transfer results to abstract Banach manifolds. Built upon the differential and integral calculus and Banach manifolds libraries in Mathlib, our work aims to lay the foundation for dynamical systems and differential geometry libraries that are general, robust, and friendly to classical mathematicians.
\end{abstract}


\section{Introduction}

The existence and uniqueness of integral curves of vector fields is a cornerstone of the theory of ordinary differential equations, with applications across dynamical systems, differential geometry, and mathematical physics. In this short paper, we present a formalisation of these results in the Lean theorem prover, working at the level of generality of Banach manifolds. Our formalisation builds on Mathlib's existing infrastructure for differential and integral calculus, Banach spaces, and smooth manifolds.

We first formalise the Picard-Lindel\"of theorem and the Gr\"onwall inequality for differential equations on Banach spaces, establishing the local existence and uniqueness of integral curves, together with corollaries such as the existence of local flows. We then transfer these results to abstract Banach manifolds, proving the local existence and uniqueness of integral curves of $C^1$ vector fields, as well as their global existence on compact manifolds. The formalisation comprises approximately 2500 lines of Lean code. To the authors' knowledge, this is the first formalisation of the existence and uniqueness of integral curves on differentiable manifolds in any proof assistant. Parts of this work have been incorporated into Mathlib.

\section{Related Work}

The formalisation in this work relies on Mathlib, a unified library of formalised mathematics in Lean~\cite{mathlib_2020}. Mathlib provides extensive libraries for differential and integral calculus on Banach spaces, including Fr\'echet derivatives and Bochner integration, as well as a library for $C^k$ differentiable manifolds and fibre bundles. Bordg and Cavalleri~\cite{Bordg_2021} describe the differential geometry library, including tangent bundles and vector fields on manifolds. Our formalisation builds directly on these libraries.

We now describe similar efforts in other proof assistants. Makarov and Spitters~\cite{Makarov_2013} formalised the Picard-Lindel\"of theorem for ODEs given by $v \colon [-a,a] \times [-K, K] \to \mathbb R$ in Coq (since renamed to Rocq), using the C-CoRN library~\cite{CruzFilipe_2004} of constructive mathematics.

The formalisation of ODEs is considerably more developed in Isabelle/HOL than in Rocq. Immler and H\"olzl~\cite{Immler_2012} formalised the Picard-Lindel\"of theorem for ODEs given by $v \colon \mathbb R \to \mathbb R^n$, together with Euler's method for the numerical approximation of a solution. Immler and Traut~\cite{Immler_2019} formalised the flow of an ODE, as well as a proof of its continuous differentiability in space via the Poincar\'e map. Immler and Tan~\cite{Immler_2020} then formalised the Poincar\'e-Bendixson theorem. Analytical convergence in these works is expressed in the theory of filters, similar to Mathlib.

We highlight some differences between our approach and those in the aforementioned works. Importantly, the dependent type theory of the Lean theorem prover~\cite{Moura2015TheLT,Moura_2021} differs from the simple type theory of Isabelle/HOL. The inability to define the tangent space at at point as a dependent type is described by Immler and Zhan~\cite{Immler_Zhan_2019} as a setback in the formalisation of differential geometry in Isabelle/HOL. In contrast, Lean's type system presents no such difficulty.

A distinguishing feature of our formalisation is its mathematical generality: by considering integral curves and flows in Banach spaces, our results apply to infinite-dimensional settings relevant to functional analysis. Previous formalisations have focused on finite-dimensional Euclidean spaces $\mathbb R^n$, often with an emphasis on computability; for instance, Immler and H\"olzl~\cite{Immler_2012} supplement their abstract formalisation with numerical methods using arbitrary-precision floats. By contrast, our work does not include any computational components, inheriting instead the philosophy of Mathlib, which strives for generality.

The unified nature of Mathlib, on which our work depends, significantly aided the formalisation process by offering libraries of disparate fields of mathematics that are coherent and interoperable. By contrast, Immler and Tan~\cite{Immler_2020} cite the combination of large and independently developed mathematical libraries as a major challenge in their formalisation. Parts of our work have already been accepted as part of Mathlib, while others are undergoing community review.

\section{Integral Curves and Flows}

We first give an overview of the definitions and statements which are the subject of this work. Let $v \colon \mathbb R \times E \to E$ be a time-dependent vector field on a Banach space $E$. An \emph{integral curve} of $v$ starting at point $x_0 \colon E$ at time $t_0 \colon \mathbb R$ is a function $\alpha \colon \mathbb R \to E$ satisfying
\begin{equation}
  \alpha(t_0) = x_0;\quad \forall t \in I,\ \partial\alpha(t) = v(t,\alpha(t)),\label{eq:integral_curve_def}
\end{equation}
where $I$ is an interval containing $t_0$. This is also called an \emph{initial value problem} for an ordinary differential equation (ODE) defined by $v$, in which case $\alpha$ is called a \emph{solution}. A (local) flow of $v$ is a one-parameter family of maps $\phi_{t} \colon E \to E$ indexed by $t \in I$, such that for all $x$ in an open subset $U \subseteq E$, the curve $t \mapsto \phi_{t}(x)$ is an integral curve of $v$ with the initial condition $\phi_0(x) = x$.

Notice that the integral curve (or flow) $\alpha$ may take arbitrary values outside its domain of $I$ (or $U \times I$), and the vector field $v$ may take arbitrary values outside $I$ and the range of $\alpha$. We will frequently exploit this ``junk value'' approach, which has the benefit of avoiding the proliferation of function types such as $v \colon I \times U \to E$, where $I$ or $U$ would change depending on the smoothness conditions on $v$.

Let $M$ be a manifold modelled on a Banach space $E$ and $v \colon \mathbb R \times M \to TM$ be a time-dependent vector field on the manifold, with $TM$ denoting the tangent bundle of $M$. An \emph{integral curve} of $v$ starting at point $x_0 \colon M$ at time $t_0 \colon \mathbb R$ is a function $\alpha \colon \mathbb R \to M$ satisfying~\eqref{eq:integral_curve_def}, where $\partial$ is now the derivation on a manifold. Similarly, a (local) flow of $v$ is a family of maps $\phi_{t} \colon M \to M$, $t \in I$, such that for all $x$ in an open subset $U \subseteq M$, the curve $t \mapsto \phi_{t}(x)$ is an integral curve of $v$ with the initial condition $\phi_0(x) = x$.

\section{Prerequisites}

The differential calculus library in Mathlib provides the predicate \texttt{HasDerivWithinAt $\alpha$ $y$ $S$ $t$}, which asserts that $\alpha$ has derivative $y$ at $t$ within the set $S$. This allows one-sided derivatives at endpoints of closed intervals. An integral curve of $v$ on $I$ is then formalised as follows:
\begin{leancode}
def IsIntegralCurveOn (α : ℝ → E) (v : ℝ → E → E) (I : Set ℝ) : Prop :=
  ∀ t ∈ I, HasDerivWithinAt α (v t (α t)) I t
\end{leancode}

For manifolds, we assume for the current work that the manifold $M$ has empty boundary, $I$ is an open interval, and $v \colon M \to TM$ is time-independent, in order to avoid issues with the boundaries of charts or directional derivatives. An integral curve on a manifold satisfies $\partial\alpha(t) = v(\alpha(t))$ for all $t \in I$, where the derivative is the manifold derivative. In Lean, this derivative is a continuous linear map from the tangent space of $\mathbb R$ at $t$ to the tangent space of $M$ at $\alpha(t)$. The vector field value $v(\alpha(t))$ is encoded as the linear map $(r \colon \mathbb R) \mapsto r \cdot v(\alpha(t))$. The formalised definition is:
\begin{leancode}
def IsMIntegralCurveOn (α : ℝ → M) (v : (x : M) → TangentSpace IM x)
    (I : Set ℝ) : Prop :=
  ∀ t ∈ I, HasMFDerivAt 𝓘(ℝ, ℝ) IM α t ((1 : ℝ →L[ℝ] ℝ).smulRight (v (α t)))
\end{leancode}
Here \texttt{$\mathscr{I}$($\mathbb{R}$, $\mathbb{R}$)} and \texttt{IM} are data structures specifying $\mathbb R$ and $M$, respectively, as real differentiable manifolds, the details of which are beyond the scope of this work.

As part of the general theory of integral curves, we provide lemmas for their translation and scaling in time, as well as variants of these definitions for local integral curves (on an arbitrarily small open time interval) and global integral curves (over all time).

For the proof of the Picard-Lindel\"of theorem, we use the integral calculus library of Mathlib. In particular, \texttt{$\int$ x in a..b, f x} denotes $\int_a^b f(x)\,dx$, formalised according to Bochner's theory of integration on Banach spaces.

All of the libraries we rely on are general enough that $E$ may be taken as a general Banach space. This is expressed in Mathlib as implicit assumptions
\begin{leancode}
variable {E : Type*} [NormedAddCommGroup E] [NormedSpace ℝ E] [CompleteSpace E]
\end{leancode}
The proof of the Picard-Lindel\"of theorem will invoke the Banach fixed point theorem, which has been formalised in Mathlib.

\section{Picard-Lindel\"of Theorem}

For a time-dependent vector field $v \colon \mathbb R \times E \to E$, the Picard-Lindel\"of theorem provides sufficient conditions for the existence and uniqueness of an integral curve $\alpha \colon \mathbb R \to E$ over a closed interval $I$ starting at the point $x_0$ at time $t_0 \in I$:
\begin{enumerate}
  \item $v(t,x)$ is $K$-Lipschitz in $x \in \overline B_a(x_0)$ for all $t \in I$, where $K \geq 0$ and $\overline B_a(x_0)$ denotes the closed ball of radius $a \geq 0$ centred at $x_0$;
  \item $v(t,x)$ is continuous in $t \in I$ for all $x \in \overline B_a(x_0)$;
  \item $\|v(t,x)\| \leq L$ for all $t \in I$ and $x \in \overline B_a(x_0)$, where $L \geq 0$;
  \item $L|t - t_0| \leq a$ for all $t \in I$.
\end{enumerate}
These conditions are explicitly stated around the point $x_0$ at which the singular integral curve starts.

To prove that there exist integral curves starting at every point $x$ in a neighbourhood $U$ of $x_0$, one could apply this theorem pointwise, choosing $x$-dependent constants $t_{\mathrm{min}}(x)$, $t_{\mathrm{max}}(x)$, $a(x)$, $K(x)$, and $L(x)$ satisfying the above conditions for each $x \in U$. However, to assemble a flow from these integral curves, one needs them to exist on a shared time interval. This requires extracting uniform constants $t_{\mathrm{min}}$, $t_{\mathrm{max}}$, $a$, $K$, and $L$ that work simultaneously for all starting points $x \in U$, making the proof obligation in downstream formalisations considerably more involved.

To remedy this, we strengthen condition 4 to
\begin{enumerate}
  \setcounter{enumi}{3}
  \item $L|t - t_0| \leq a - r$ for all $t \in I$, where $r \geq 0$.
\end{enumerate}
This guarantees the existence of a local flow of $v$ on $\overline B_r(x_0)$. Now, downstream proofs using the Picard-Lindel\"of theorem for the existence of flows on a neighbourhood of $x_0$ only has to provide one set of constants $t_{\mathrm{min}}$, $t_{\mathrm{max}}$, $a$, $r$, $K$, and $L$ and show four conditions around a fixed base point $x_0$, rather than a whole family of conditions indexed by starting points near $x_0$. To recover the usual statement for a single integral curve starting at $x_0$, one simply sets $r = 0$.

Since the conditions on $v$ are four long statements, we formalise them as a proposition-valued structure, so that each condition may be referred to by name:
\begin{leancode}
structure IsPicardLindelof (v : ℝ → E → E) {tmin tmax : ℝ} (t₀ : Icc tmin tmax)
    (x₀ : E) (a r L K : ℝ≥0) : Prop where
  lipschitzOnWith : ∀ t ∈ Icc tmin tmax, LipschitzOnWith K (v t) (closedBall x₀ a)
  continuousOn : ∀ x ∈ closedBall x₀ a, ContinuousOn (v · x) (Icc tmin tmax)
  norm_le : ∀ t ∈ Icc tmin tmax, ∀ x ∈ closedBall x₀ a, ‖v t x‖ ≤ L
  mul_max_le : L * max (tmax - t₀) (t₀ - tmin) ≤ a - r
\end{leancode}
where \texttt{Icc tmin tmax} denotes the closed interval $[t_{\mathrm{min}}, t_{\mathrm{max}}]$ The statement of the existence theorem is then
\begin{leancode}
theorem exists_eq_isIntegralCurveOn
    (hf : IsPicardLindelof f t₀ x₀ a r L K) (hx : x ∈ closedBall x₀ r) :
    ∃ α : ℝ → E, α t₀ = x ∧ IsIntegralCurveOn α f (Icc tmin tmax) := ...
\end{leancode}

Our formalisation of the proof of the Picard-Lindel\"of theorem follows the standard procedure of defining a complete space of functions on which the Picard iteration acts as a contracting map, and then extracting the integral curve as the limit of the iteration using the Banach fixed point theorem. The space of functions may be taken to be the space of bounded continuous functions $\mathbb R \to \mathbb R^n$ as by Immler and H\"olzl~\cite{Immler_2012}, or the space of $\mu$-uniformly continuous functions as by Makarov and Spitters~\cite{Makarov_2013}. but we have chosen to use the space of $L$-Lipschitz functions starting within distance $r$ of the point $x_0$, where $r$ and $L$ are the same as the ones in \texttt{IsPicardLindelof}:
\begin{leancode}
structure FunSpace (t₀ : Icc tmin tmax) (x₀ : E) (r L : ℝ≥0) where
  toFun : Icc tmin tmax → E
  lipschitzWith : LipschitzWith L toFun
  mem_closedBall₀ : toFun t₀ ∈ closedBall x₀ r
\end{leancode}
There are many possible and mathematically inequivalent choices for this function space, with no consequence on the generality of the result. In practice, the choice depends on how conveniently it can be expressed using existing libraries, and how easily one can verify the conditions required by the Banach fixed point theorem: that the space is complete and that the Picard iteration maps the space into itself.

The condition \lean{mem_closedBall₀} requires only that the curve starts within distance $r$ of a fixed base point $x_0$. This allows Picard iterations with different starting points $x \in \overline B_r(x_0)$ to operate on the same \lean{FunSpace}, simplifying the formalisation of the existence of the flow.

The Banach fixed point theorem requires that \texttt{FunSpace} be a complete metric space under the supremum norm. Therefore, we explicitly specify the domain as $I \to E$ rather than use the junk value approach, which would have produced a pseudometric space instead.

The Picard iteration is formalised as an operator on \texttt{FunSpace}. First, the function \texttt{picard} takes a curve $\alpha \colon \mathbb R \to E$ (assuming junk values outside of $I$) and returns the next iterate $t \mapsto x_0 + \int_{t_0}^t v(\tau, \alpha(\tau))\,d\tau$:
\begin{leancode}
noncomputable def picard (f : ℝ → E → E) (t₀ : ℝ) (x₀ : E) (α : ℝ → E) :
    ℝ → E := fun t ↦ x₀ + ∫ τ in t₀..t, f τ (α τ)
\end{leancode}
Second, we lift \texttt{picard} to \texttt{next}, which operates on \texttt{FunSpace} and produces another element of \texttt{FunSpace}. This is done with the help of \texttt{compProj}, which extends a function on $I$ to all of $\mathbb R$ by mapping points outside $I$ to the nearest endpoint.
\begin{leancode}
noncomputable def next (hf : IsPicardLindelof f t₀ x₀ a r L K)
    (hx : x ∈ closedBall x₀ r) (α : FunSpace t₀ x₀ r L) : FunSpace t₀ x₀ r L where
  toFun t := picard f t₀ x α.compProj t
  lipschitzWith := ...
  mem_closedBall₀ := ...
\end{leancode}

The remainder of the proof shows that \texttt{next}, when iterated a number of times, is a contracting map on \texttt{FunSpace}, and the Banach fixed point theorem is then applied to obtain a fixed point of \texttt{next}. The existence of an integral curve follows.

A local flow $E \times \mathbb R \to E$ is then constructed from the existence of a collection of integral curves $\mathbb R \to E$ with different starting points (invoking the axiom of choice):
\begin{leancode}
theorem exists_forall_mem_closedBall_eq_isIntegralCurveOn
    (hf : IsPicardLindelof f t₀ x₀ a r L K) :
    ∃ α : E → ℝ → E, ∀ x ∈ closedBall x₀ r, α x t₀ = x ∧
      IsIntegralCurveOn (α x) f (Icc tmin tmax) := ...
\end{leancode}

We also prove a number of corollaries. The flow is continuous on $I \times U$. If the vector field $v \colon \mathbb R \times E \to E$ is $C^k$ on $I \times U$, then its integral curve is $C^k$ in time. If a time-independent vector field $v \colon E \to E$ is $C^1$ at $x_0$, then it satisfies the conditions of the Picard-Lindel\"of theorem at $x_0$ for some values of $a,L,K$.

\section{Gr\"onwall Inequality}

The uniqueness of integral curves in Banach spaces is established through the Gr\"onwall inequality, which places an exponential upper bound on the distance between two $\epsilon$-approximate integral curves of the same vector field starting a distance $\delta$ apart:
\begin{leancode}
noncomputable def gronwallBound (δ K ε x : ℝ) : ℝ :=
  if K = 0 then δ + ε * x else δ * exp (K * x) + ε / K * (exp (K * x) - 1)

theorem dist_le_of_approx_trajectories_ODE_of_mem
    (hv : ∀ t ∈ Ico a b, LipschitzOnWith K (v t) (s t))
    (hf : ContinuousOn f (Icc a b)) (hfs : ∀ t ∈ Ico a b, f t ∈ s t)
    (hf' : ∀ t ∈ Ico a b, HasDerivWithinAt f (f' t) (Ici t) t)
    (f_bound : ∀ t ∈ Ico a b, dist (f' t) (v t (f t)) ≤ εf)
    ... -- similar assumptions on g
    (ha : dist (f a) (g a) ≤ δ) :
    ∀ t ∈ Icc a b, dist (f t) (g t) ≤ gronwallBound δ K (εf + εg) (t - a) := ...
\end{leancode}
Our formalisation, which follows Hubbard and West~\cite[Section 4.5]{hubbard_2013_differential}, relies on an existing formalised proof in Mathlib of a bound on a continuous function with an estimate of its derivative. The uniqueness theorem then follows as a corollary.

This approach differs from the Isabelle/HOL formalisations~\cite{Immler_2012,Immler_2019}, where uniqueness is obtained as part of the Picard-Lindel\"of construction via the Banach fixed point theorem. By contrast, Mathlib's Gr\"onwall-based proof of uniqueness does not depend on integral operators or fixed point theorems. This separation is advantageous for modularity: once existence of a solution is established by any means, uniqueness can be proven using only elementary estimates, without any dependence on Banach fixed point and integral libraries.

For ease of use, we then provide a number of corollaries where the assumptions take different forms and where integral curves extend ``backwards'' before the initial time.

\section{Transferring to Manifolds}

Given a vector field on a manifold, $v \colon M \to TM$, that is $C^1$ at $x_0 \colon M$, we now wish to show that there exists an integral curve $\alpha \colon \mathbb R \to M$ of $v$ starting at $x_0$ at time $t_0$, over some open interval $I = (t_0 - \epsilon, t_0 + \epsilon)$, where $\epsilon > 0$.

In Mathlib, the tangent space \texttt{TangentSpace IM x} at each point $x \colon M$ is defined as a type synonym of the model space $E$. A vector field is a section \texttt{v : (x : M) $\to$ TangentSpace IM x}, where \texttt{v x} is the value of the vector when expressed via the local trivialisation at $x$.

We now schematically describe our formalisation using more familiar mathematical notation for legibility. For each $x \colon M$, let $\phi_x \colon M \to E$ denote the (preferred) local chart at $x$, and let $\phi^{TM}_x \colon TM \to E \times E$ denote the local trivialisation of $TM$ at $x$. Note again the junk value approach, where only the values within the base set $U_x$ of the local chart at $x$ matter. The local trivialisations of $TM$ are induced by the local charts of $M$, and they share the same base sets $U_x$. For each $x \colon M$ and $x' \in U_x$,
$$[\phi^{TM}_x(v(x'))]_1 = \phi_x(x'),\quad [\phi^{TM}_x(v(x'))]_2 = \partial(\phi_x \circ \phi_{x'}^{-1}) [\phi^{TM}_{x'}(v(x'))]_2,$$
where $[\phi^{TM}_{x'}(v(x'))]_2$ is exactly the value of \texttt{v x' : TangentSpace IM x'}.

A vector field $v \colon M \to TM$ is said to be $C^1$ at $x$ if $\phi^{TM}_x \circ v \circ \phi_x^{-1} \colon E \to E\times E$ is $C^1$ at $\phi_x(x)$. In particular, $z \mapsto [\phi^{TM}_x \circ v \circ \phi_x^{-1}(z)]_2$ is a vector field $E \to E$ on a Banach space that is $C^1$ at $\phi_x(x)$. For any starting time $t_0$, we can then apply a corollary of the Picard-Lindel\"of theorem to obtain an integral curve $\tilde\alpha \colon \mathbb R \to E$ such that
$$\tilde\alpha(t_0) = \phi_x(x);\quad \exists \epsilon > 0,\ \forall t \in (t_0 - \epsilon, t_0 + \epsilon),\ \partial\tilde\alpha(t) = [\phi^{TM}_x \circ v \circ \phi_x^{-1}(\tilde\alpha(t))]_2.$$
The integral curve on the manifold is then constructed as $\alpha := \phi_x^{-1} \circ \tilde\alpha$, with $\epsilon$ appropriately shrunken to ensure $\tilde\alpha(t) \in \phi_x(U_x)$.

It remains to show that the manifold derivative of $\alpha$ coincides with $v$:
$$\partial_t(\phi_{\alpha(t)} \circ \alpha(t)) = [v(\alpha(t))]_2.$$
The crucial step is to multiply the right hand side by the identity, $[\partial (\phi_x \circ \phi_{\alpha(t)}^{-1})]^{-1} \partial (\phi_x \circ \phi_{\alpha(t)}^{-1})$, before applying the chain rule on the left hand side.

The proof that integral curves on manifolds are unique proceeds similarly. To prove that they are unique on any open interval $I \subset \mathbb R$, we follow Lee~\cite[Theorem 9.12(a)]{lee_2012_introduction} by showing that the set $\{t \mid \alpha(t) = \alpha'(t)\}$ is clopen in $I$ for any two integral curves $\alpha,\alpha'$ of $v$ with the same starting point. The global uniqueness of integral curves then follows.

As a corollary, we also formalise the proof that $C^1$ vector fields on compact manifolds admit global integral curves.

\section{Ongoing and Future Work}

The current work is part of a general effort to build up formalised libraries of dynamical systems and differential geometry, which has so far received less attention and contribution than, for instance, libraries of algebra and category theory.

We are in the process of formalising the theorem that the flow of a $C^k$ vector field on a Banach space $E$ is $C^k$ on its domain (in $E$). The version for $C^1$ vector fields where $E = \mathbb R^n$ has already been formalised by Immler and Traut~\cite{Immler_2019} using the Poincar\'e map, but we intend to formalise a more general result. We plan to follow the proof of Lang~\cite[page 81 onwards]{lang_2012_fundamentals} based on the implicit function theorem on Banach spaces, since it can be adapted to the case of $H^k$ Sobolev spaces, which have not yet been formalised in Lean.

Since Mathlib already contains libraries of Lie algebras and Lie groups, this work enables the formalisation of the exponential map and hence the correspondence between representations of Lie algebras and Lie groups.

Although Riemannian geometry has not been formalised in Mathlib, this work is essential for stating the existence of geodesics on Riemannian manifolds.

The case where the integral curve may venture to the boundary of the manifold requires special treatment, for example, of the direction of the vector field on the boundary. We will incorporate such cases where we can and postpone ones that need the invariance of the boundary of a manifold, which has not been formalised in Lean.

This work demonstrates that topics in differential geometry are amenable to formalisation, without any concerns about constructivism or computability. The authors hope that this will inspire more interest among mathematicians in the formalisation of differential geometry, with the aim of reaching major results and cutting-edge research.



\bibliography{references}

@InProceedings{Immler_2012,
author="Immler, Fabian
and H{\"o}lzl, Johannes",
editor="Beringer, Lennart
and Felty, Amy",
title="Numerical Analysis of Ordinary Differential Equations in {Isabelle/HOL}",
booktitle="Interactive Theorem Proving",
year="2012",
publisher="Springer Berlin Heidelberg",
address="Berlin, Heidelberg",
pages="377--392",
isbn="978-3-642-32347-8"
}

@InProceedings{Makarov_2013,
author="Makarov, Evgeny
and Spitters, Bas",
editor="Blazy, Sandrine
and Paulin-Mohring, Christine
and Pichardie, David",
title="The {Picard} Algorithm for Ordinary Differential Equations in {Coq}",
booktitle="Interactive Theorem Proving",
year="2013",
publisher="Springer Berlin Heidelberg",
address="Berlin, Heidelberg",
pages="463--468",
isbn="978-3-642-39634-2"
}

@article{Immler_2019,
author = {Immler, Fabian and Traut, Christoph},
title = {The Flow of {ODEs}: Formalization of Variational Equation and {Poincar\'{e}} Map},
year = {2019},
issue_date = {February  2019},
publisher = {Springer-Verlag},
address = {Berlin, Heidelberg},
volume = {62},
number = {2},
issn = {0168-7433},
url = {https://doi.org/10.1007/s10817-018-9449-5},
doi = {10.1007/s10817-018-9449-5},
journal = {J. Autom. Reason.},
month = feb,
pages = {215–236},
numpages = {22}
}

@inproceedings{Moura2015TheLT,
  title={The {Lean} Theorem Prover (System Description)},
  author={Leonardo Mendonça de Moura and Soonho Kong and Jeremy Avigad and Floris van Doorn and Jakob von Raumer},
  booktitle={CADE},
  year={2015},
  url={https://api.semanticscholar.org/CorpusID:232990}
}

@InProceedings{CruzFilipe_2004,
author="Cruz-Filipe, Lu{\'i}s
and Geuvers, Herman
and Wiedijk, Freek",
editor="Asperti, Andrea
and Bancerek, Grzegorz
and Trybulec, Andrzej",
title="{C-CoRN, the Constructive Coq Repository at Nijmegen}",
booktitle="Mathematical Knowledge Management",
year="2004",
publisher="Springer Berlin Heidelberg",
address="Berlin, Heidelberg",
pages="88--103",
isbn="978-3-540-27818-4"
}

@inproceedings{Immler_2020,
author = {Immler, Fabian and Tan, Yong Kiam},
title = {The {Poincar\'{e}-Bendixson} theorem in {Isabelle/HOL}},
year = {2020},
isbn = {9781450370974},
publisher = {Association for Computing Machinery},
address = {New York, NY, USA},
url = {https://doi.org/10.1145/3372885.3373833},
doi = {10.1145/3372885.3373833},
booktitle = {Proceedings of the 9th ACM SIGPLAN International Conference on Certified Programs and Proofs},
pages = {338–352},
numpages = {15},
location = {New Orleans, LA, USA},
series = {CPP 2020}
}

@inproceedings{Immler_Zhan_2019,
author = {Immler, Fabian and Zhan, Bohua},
title = {Smooth manifolds and types to sets for linear algebra in {Isabelle/HOL}},
year = {2019},
isbn = {9781450362221},
publisher = {Association for Computing Machinery},
address = {New York, NY, USA},
url = {https://doi.org/10.1145/3293880.3294093},
doi = {10.1145/3293880.3294093},
booktitle = {Proceedings of the 8th ACM SIGPLAN International Conference on Certified Programs and Proofs},
pages = {65–77},
numpages = {13},
location = {Cascais, Portugal},
series = {CPP 2019}
}

@book{lee_2012_introduction,
  author = {Lee, John M},
  publisher = {Springer},
  title = {Introduction to Smooth Manifolds},
  year = {2012}
}

@book{hubbard_2013_differential,
  author = {Hubbard, John H and West, Beverly H},
  month = {11},
  publisher = {Springer},
  title = {Differential Equations: A Dynamical Systems Approach},
  year = {2013}
}

@book{lang_2012_fundamentals,
  author = {Lang, Serge},
  month = {12},
  publisher = {Springer},
  title = {Fundamentals of Differential Geometry},
  year = {2012}
}

@inproceedings{mathlib_2020,
author = {The mathlib Community},
title = {The {Lean} mathematical library},
year = {2020},
isbn = {9781450370974},
publisher = {Association for Computing Machinery},
address = {New York, NY, USA},
url = {https://doi.org/10.1145/3372885.3373824},
doi = {10.1145/3372885.3373824},
booktitle = {Proceedings of the 9th ACM SIGPLAN International Conference on Certified Programs and Proofs},
pages = {367–381},
numpages = {15},
location = {New Orleans, LA, USA},
series = {CPP 2020}
}

@InProceedings{Moura_2021,
author="Moura, Leonardo de
and Ullrich, Sebastian",
editor="Platzer, Andr{\'e}
and Sutcliffe, Geoff",
title="The Lean 4 Theorem Prover and Programming Language",
booktitle="Automated Deduction -- CADE 28",
year="2021",
publisher="Springer International Publishing",
address="Cham",
pages="625--635",
isbn="978-3-030-79876-5"
}

@article{Bordg_2021,
  author       = {Anthony Bordg and
                  Nicol{\`{o}} Cavalleri},
  title        = {Elements of Differential Geometry in {Lean}: {A} Report for Mathematicians},
  journal      = {CoRR},
  volume       = {abs/2108.00484},
  year         = {2021},
  url          = {https://arxiv.org/abs/2108.00484},
  eprinttype    = {arXiv},
  eprint       = {2108.00484},
  bibsource    = {dblp computer science bibliography, https://dblp.org}
}

\end{document}